\documentclass[12pt,doublespacing]{article}

\usepackage{fullpage}
\usepackage{setspace}
\usepackage{authblk}

\usepackage{cite}
\usepackage{hyperref} 

\usepackage{graphicx}

\usepackage{amsmath,amssymb,amsfonts,physics,esint,mathtools}
\usepackage{algorithmic,algorithm}
\usepackage{xcolor}

\usepackage{float}
\usepackage{subfig}
\usepackage{makecell}
\usepackage{rotating}
\usepackage{graphicx}
\usepackage{textcomp}

\usepackage{enumitem}

\def\Xint#1{\mathchoice
	{\XXint\displaystyle\textstyle{#1}}%
	{\XXint\textstyle\scriptstyle{#1}}%
	{\XXint\scriptstyle\scriptscriptstyle{#1}}%
	{\XXint\scriptscriptstyle\scriptscriptstyle{#1}}%
	\!\int}
\def\XXint#1#2#3{{\setbox0=\hbox{$#1{#2#3}{\int}$}
		\vcenter{\hbox{$#2#3$}}\kern-.5\wd0}}

\def\dashint{\Xint-}
\def\BibTeX{{\rm B\kern-.05em{\sc i\kern-.025em b}\kern-.08em
    T\kern-.1667em\lower.7ex\hbox{E}\kern-.125emX}}
\begin{document}
%

\title{A Multi-Frequency Iterative Method for Reconstruction of Rough Surfaces Separating Two Penetrable Media} 

\author[1,2]{Ahmet Sefer}
\author[3]{Ali Yapar}
\author[1]{Hakan Bagci\vspace{0.5cm}}

\affil[1]{Electrical and Computer Engineering (ECE) Program 
\authorcr Computer, Electrical, and Mathematical Science and Engineering (CEMSE) Division, 
\authorcr King Abdullah University of Science and Technology (KAUST) 
\authorcr Thuwal 23955, Saudi Arabia
\authorcr e-mail: $\lbrace$ahmet.sefer, hakan.bagci$\rbrace$@kaust.edu.sa\vspace{0.5cm}}

\affil[2]{Electrical and Electronics Engineering 
\authorcr Fevziye Schools Foundation Isik University, Istanbul 34882, Turkey
\authorcr e-mail: ahmet.sefer@isikun.edu.tr\vspace{0.5cm}}

\affil[3]{Electronics and Communications Engineering 
\authorcr Istanbul Technical University, Istanbul 34469, Turkey
\authorcr e-mail: yapara@itu.edu.tr\vspace{0.5cm}}


\date{}
\maketitle
\newpage

\begin{abstract}
A numerical scheme that uses multi-frequency Newton iterations to reconstruct a rough surface profile between two dielectric media is proposed. At each frequency sample, the scheme employs Newton iterations to solve the nonlinear inverse scattering problem. At every iteration, the Newton step is computed by solving a linear system that involves the Frechet derivative of the integral operator, which represents the scattered fields, and the difference between these fields and the measurements. This linear system is regularized using the Tikhonov method. The multi-frequency data is accounted for in a recursive manner. More specifically, the profile reconstructed at a given frequency is used as an initial guess for the iterations at the next frequency. The effectiveness of the proposed method is validated through numerical examples, which demonstrate its ability to accurately reconstruct surface profiles even in the presence of measurement noise. The results also show the superiority of the multi-frequency approach over single-frequency reconstructions, particularly in terms of handling surfaces with sharp variations.

\end{abstract}

\newpage
\section{Introduction}\label{sec:introduction}
Reconstruction of inaccessible rough surfaces from measured scattered electromagnetic fields is a subject of significant interest in various engineering disciplines, such as remote sensing~\cite{RemoteSensing1,RemoteSensing2,RemoteSensing3,RemoteSensing4_optical}, optical system measurement~\cite{OpticalSystMea}, subsurface imaging~\cite{SubSurfaceImaging1,SubsurfaceImaging2,multiTGRS2015}, ultrasonic applications like wall-thickness measurement~\cite{UltrasonicWallThickness}, damage detection~\cite{UltrasonicDamageDetection,UltrasonicLaserImaging}, and nondestructive testing~\cite{UltrasonicNDT,UltrasonicNDT2,UltrasonicsNDT3}. This reconstruction requires solving an inverse problem where the scattered fields are represented as convolutions of the Green functions of the background media with the fields on the unknown surface profile~\cite{tsang}. This inverse problem is inherently ill-posed due to the contamination of the measured scattered field data by noise and the ``smoothing'' effect introduced by the convolution integrals~\cite{colton,kressReview2018}. Furthermore, the scattered fields are nonlinear functions of the unknown surface profile~\cite{colton,kressReview2018}. The ill-posedness and nonlinearity of the inverse problem make the reconstruction of the surface profile a highly challenging task. 

Among the methods that are developed to address these challenges, semi-analytical approaches that rely on Kirchhoff~\cite{dolcetti2020}, small-perturbation~\cite{wombellSPM} and Rytov~\cite{Rytov} approximations, or low-order expansion of fields~\cite{bao2} and fully numerical approaches that rely on reverse time migration (RTM)~\cite{asRTMTGRS} can be considered ``direct'' solution techniques, i.e., they are non-iterative. 

The other group of solution techniques~\cite{chen, chen2,Spivack2024, Bao13Landweber, Bao2016,chorfi2011,asTGRS,Li2015,asTGRS2,asGRSL} minimize the error between the measured scattered fields and the scattered fields of a predicted profile that is updated iteratively to ``linearize'' the inverse problem. Often, the regularization is applied at every iteration to alleviate the ill-posedness. The method described in~\cite{chen,chen2,Spivack2024} iteratively updates the derivative of the field on the surface and the surface profile that are coupled via the convolution integral and a simple relationship that is obtained under the assumption of grazing incident field. In~\cite{Bao13Landweber} and~\cite{Bao2016}, Landweber iterations are used for linearization and regularization of the reconstruction of periodic gratings and rough surfaces from phaseless data, respectively. In contrast, the method proposed in~\cite{chorfi2011} inverts the full scattered fields (with phase and amplitude) to reconstruct rough surfaces separating two dielectric media. A Newton method is used for linearization while the regularization at every Newton iteration is carried out by applying the truncated conjugate gradient method to the normal equation of the Newton update.

The other group of solution techniques~\cite{chen, chen2,Spivack2024, Bao13Landweber, Bao2016,chorfi2011,asTGRS,Li2015,asTGRS2,asGRSL} minimize the error between the measured scattered fields and the scattered fields of a predicted profile that is updated iteratively to ``linearize'' the inverse problem. Often, the regularization is applied at every iteration to alleviate the ill-posedness. The method described in~\cite{chen,chen2,Spivack2024} iteratively updates the derivative of the field on the surface and the surface profile that are coupled via the convolution integral and a simple relationship that is obtained under the assumption of grazing incident field. In~\cite{Bao13Landweber} and~\cite{Bao2016}, Landweber iterations are used for linearization and regularization of the reconstruction of periodic gratings and rough surfaces from phaseless data, respectively. In contrast, the method proposed in~\cite{chorfi2011} inverts the full scattered fields (with phase and amplitude) to reconstruct rough surfaces separating two dielectric media. A Newton method is used for linearization while the regularization at every Newton iteration is carried out by applying the truncated conjugate gradient method to the normal equation of the Newton update.

Similarly, in~\cite{asTGRS}, a Newton method is used for the reconstruction of rough surfaces separating two dielectric media from full scattered-field measurements. At every iteration, a linear system in the Newton step of the unknown profile is constructed. This linear system involves the Frechet derivative of the convolution operator, which is used in the representation of the fields scattered from the profile updated at that iteration and the difference between these scattered fields and the measurements. To alleviate the ill-posedness of this linear system, Tikhonov regularization is applied before it is solved for the Newton step. In~\cite{Li2015}, a similar iterative method is developed for the reconstruction of sound-soft rough surfaces of acoustics. In~\cite{asTGRS2} and~\cite{asGRSL}, the Newton method in~\cite{asTGRS} is extended for reconstruction using phaseless data and reconstruction of impedance surfaces, respectively.

The performance of these iterative approaches can be improved using multi-frequency/multi-resolution based techniques since the use of multi-frequency/scale data alleviates the effects of ill-conditioning, reduces the occurrence of false solutions, and helps to avoid local minima of the minimization problem by mitigating the effects of non-linearity~\cite{multiExplain, multiFreqHop, multiTGRSmultiHope}. Indeed, these techniques have been used to improve solutions of inverse scattering problems in a range applications changing from ground penetrating radar (GPR)~\cite{multiTGRS2015,multiJimagingGPR} to microwave imaging~\cite{ multiMicrowaveTrans}, non-destructive testing~\cite{multiNDT} and diffraction tomography~\cite{multiTGRS2022}. In~\cite{multiTGRS2015}, a multi-scale and multi-frequency approach is used to iteratively reconstruct the scatterer profile from time-domain GPR data. In~\cite{multiJimagingGPR}, a multi-frequency contrast source imaging (CSI) method that exploits multi-view wide band GPR data is developed to reconstruct pixel-sparse subsurface objects. In~\cite{multiMicrowaveTrans}, a recursive multi-scale approach is used in conjunction with the contradiction integral equation to retrieve the unknown relative permittivity of a complex-shaped strong scatterer.

In this work, the Tikhonov-regularized Newton iterative scheme, which is proposed in~\cite{asTGRS} to reconstruct a rough surface separating two dielectric media from full scattered-field measurements, is extended to account for multi-frequency data. Execution over multiple frequency samples allows for this method to capture more details about the roughness of the surface as the wavelength gets smaller at each successive frequency. In addition to higher resolution, multi-frequency execution increases the robustness of the algorithm since the Newton iterations at a given frequency use the reconstruction at the previous frequency as an initial guess (which leads to increased convergence). 

The rest of the paper is organized as follows: Section~\ref{sec:formulation} expounds on the formulation underlying the proposed multi-frequency Newton iterations. Section~\ref{sec:setup} describes the set-up of the problem. Section~\ref{sec:scattering_problem} provides the formulation for the forward scattering problem in terms of integral equations. Section~\ref{sec:inv_scattering_problem} describes the linearization of the nonlinear inverse scattering problem via Newton iterations and its regularization via the Tikhonov method. Section~\ref{sec:multi-freq} provides the final form of the proposed method in the form of an algorithm and describes how multi-frequency data is accounted for. Comprehensive numerical results are provided in Section~\ref{sec:numerical_results} to demonstrate the effects of the simulation and problem parameters on the reconstruction accuracy. Finally, conclusions follow in Section~\ref{sec5conc}.

\section{Formulation}\label{sec:formulation}
\subsection{Problem Setup}\label{sec:setup}
\begin{figure}[ht!]
	\centering
	\includegraphics[width=0.8\columnwidth]{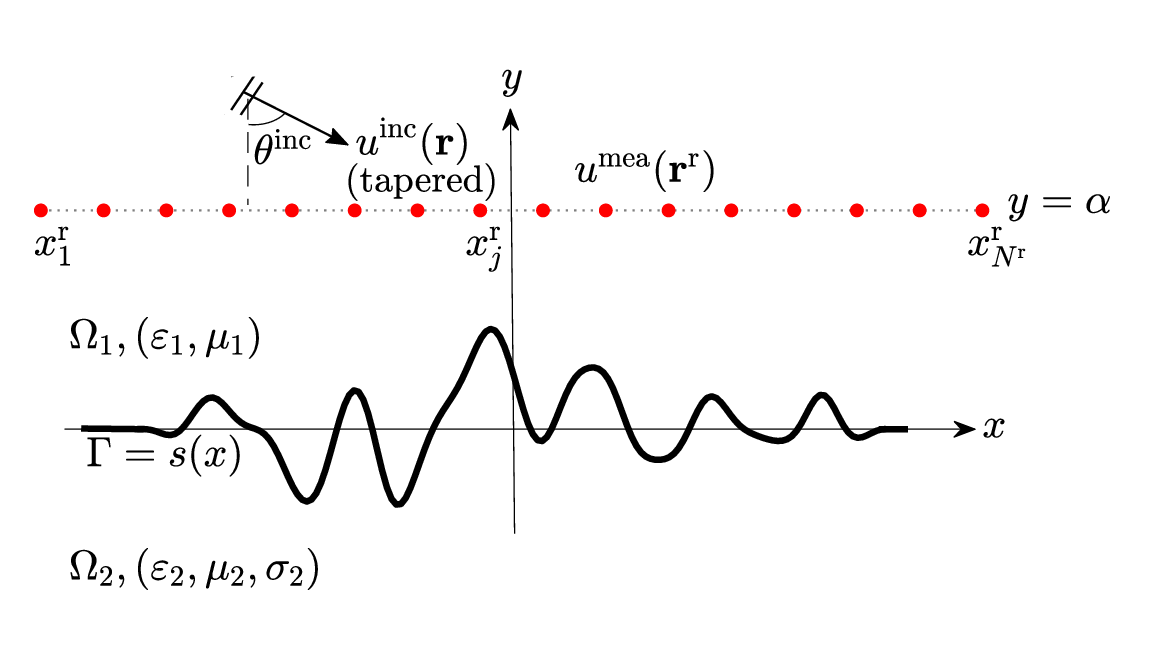}
	\caption{2D scattering problem involving a rough surface separating two dielectric media.} \label{fig:figure1}
\end{figure}

Fig.~\ref{fig:figure1} describes the two-dimensional (2D) scattering problem involving a rough surface that separates two penetrable media. It is assumed that $\Omega_1$ is lossless and $\Omega_2$ is lossy with finite conductivity. The permeability, the permittivity, and the wavenumber in $\Omega_1$ are denoted by $\mu_1$, $\varepsilon_1$, and $k_1$ and the permeability, the permittivity, the conductivity, and the wavenumber in $\Omega_2$ are denoted by  $\mu_2$, $\varepsilon_2$, $\sigma_2$, and $k_2$, respectively. The rough surface separating $\Omega_1$ and $\Omega_2$ is denoted by $\Gamma$ and expressed using a continuous height function $y = s(x)$, $L/2\geq x\geq-L/2$. Since $\Gamma$ is assumed to be of finite length, a ``traditional'' plane wave excitation gives rise to diffraction on the edges of the surface. Therefore, a plane wave with the Thorsos taper is used as excitation~\cite{thorsos}. Assuming that the plane wave originates in $\Omega_1$, the incident field is expressed using
\begin{equation}
u^{\mathrm{inc}}(\mathbf{r})=e^{\mathrm{i}k_1\hat{\mathbf{k}}^{\mathrm{inc}} \cdot \mathbf{r}} e^{-\left(\frac{x+y\tan \theta^{\mathrm{inc}}}{g}\right)^2} e^{\mathrm{i}\left(\xi(\mathbf{r}) k_1\hat{\mathbf{k}}^{\mathrm{inc}}\cdot \mathbf{r}\right)}.
\label{Tapered}
\end{equation}
In~\eqref{Tapered}, $\mathbf{r}=(x,y)$ is the location vector in the 2D space, $\hat{\mathbf{k}}^{\mathrm{inc}} = (\sin\theta^{\mathrm{inc}},-\cos\theta^{\mathrm{inc}})$ is the direction of propagation, $\theta^{\mathrm{inc}}$ is the angle of incidence, and the second and the third exponents are the decay factor and the correction term associated with the Thorsos taper, respectively. The decay factor is defined such that $u^{\mathrm{inc}}(\mathbf{r})$ decays in the direction perpendicular to $\hat{\mathbf{k}}^{\mathrm{inc}}$. The correction term, where the function $\xi(\mathbf{r})$ is defined as
\begin{equation}
\xi(\mathbf{r})=\left[\left(\frac{2\left(x+y\tan \theta^{\mathrm{inc}}\right)^2}{g^2}\right)-1\right] \frac{1}{\left(k_1 g \cos \theta^{\mathrm{inc}}\right)^2}
\label{taper_function}
\end{equation}
ensures that $u^{\mathrm{inc}}(\mathbf{r})$ satisfies the scalar Helmholtz equation to order $1 /(k_1 g\cos\theta^{\mathrm{inc}})^3$~\cite{tsang}. In~\eqref{Tapered} and~\eqref{taper_function}, the parameter $g$ controls the width of the taper. 

\subsection{Forward Scattering Problem}\label{sec:scattering_problem}
Let $u_1(\mathbf{r})$ and $u_2(\mathbf{r})$ represent the total field in $\Omega_1$ and $\Omega_2$, respectively. Using equivalence and extinction theorems~\cite{ExtinctionThrm}, one can obtain the integral representations of $u_1(\mathbf{r})$ and $u_2(\mathbf{r})$ as
\begin{subequations}
\begin{align}
 u_1(\mathbf{r})=&u^{\mathrm{inc}}(\mathbf{r})+\int_{\Gamma}\Big[ K_1(\mathbf{r},\mathbf{r}^\prime)u_1(\mathbf{r}^\prime)-G_1(\mathbf{r},\mathbf{r}^\prime)v_1(\mathbf{r}^\prime)\Big]dl'
\mathbf{r}\in \Omega_1\label{SIEa}\\ 
 u_2(\mathbf{r})=&-\int_{\Gamma}\Big[K_2(\mathbf{r},\mathbf{r}^\prime)u_2(\mathbf{r}^\prime)-G_2(\mathbf{r},\mathbf{r}^\prime)v_2(\mathbf{r}^\prime)\Big]dl'
\mathbf{r}\in \Omega_2.\label{SIEb} 
\end{align}
\label{SIE}
\end{subequations}
\noindent Here, $G_m(\mathbf{r},\mathbf{r}^\prime)=(\mathrm{i}/4)H_0^{(1)}(k_m\abs{\mathbf{r}-\mathbf{r}^\prime})$, $m\in\{1,2\}$ is the fundamental solution (Green function) of the scalar Helmholtz equation in 2D unbounded space with wavenumber $k_m$, $H_0^{(1)}(.)$ is the Hankel function of the first kind and order zero,  and $K_m(\mathbf{r},\mathbf{r}^\prime)=\hat{\mathbf{n}}(\mathbf{r}^\prime)\cdot \nabla'G_m(\mathbf{r},\mathbf{r}^\prime)$, $\mathbf{r}^\prime\in \Gamma$ is the derivative of $G_m(\mathbf{r},\mathbf{r}^\prime)$ with respect to surface unit normal vector $\hat{\mathbf{n}}(\mathbf{r}^\prime)$. Note that $\hat{\mathbf{n}}(\mathbf{r})$, $\mathbf{r}\in \Gamma$ points from $\Omega_2$ to $\Omega_1$. Similarly, $v_m(\mathbf{r})=\hat{\mathbf{n}}(\mathbf{r})\cdot\nabla u_m(\mathbf{r})$, $\mathbf{r} \in \Gamma$ is the derivative of $u_m(\mathbf{r})$ with respect to $\hat{\mathbf{n}}(\mathbf{r})$. The fields $u_1(\mathbf{r})$ and $u_2(\mathbf{r})$ and their normal derivatives $v_1(\mathbf{r})$ and $v_2(\mathbf{r})$ are continuous on $\Gamma$:
\begin{subequations}
\begin{align}
\label{EqBCsa}    u_1(\mathbf{r})&=u_2(\mathbf{r})=u(\mathbf{r}),\, \mathbf{r}\in\Gamma \\
\label{EqBCsb} v_1(\mathbf{r})&=v_2(\mathbf{r})=v(\mathbf{r}), \, \mathbf{r}\in\Gamma.
\end{align}
\label{EqBCs}
\end{subequations}
Inserting~\eqref{EqBCsa} and~\eqref{EqBCsb} into~\eqref{SIEa} and~\eqref{SIEb} and letting $\mathbf{r}\to \Gamma$ yield a coupled system of two integral equations: 
\begin{subequations}
\begin{align}
\frac{1}{2}u(\mathbf{r})-\dashint_\Gamma K_1(\mathbf{r},\mathbf{r}')u(\mathbf{r'})\,dl'
&+\int_\Gamma G_1(\mathbf{r},\mathbf{r}')v(\mathbf{r'})\, dl'
=u^{\mathrm{inc}}(\mathbf{r}),\,\mathbf{r}\in \Gamma
\label{SIE2a} \\
\frac{1}{2}u(\mathbf{r})+\dashint_\Gamma K_2(\mathbf{r},\mathbf{r}')u(\mathbf{r'})\,dl'&-\int_\Gamma G_2(\mathbf{r},\mathbf{r}')v(\mathbf{r'}) dl'
=0,\,\mathbf{r}\in \Gamma \label{SIE2b}
\end{align}
\label{SIE2}
\end{subequations}
Symbol ``$-$'' shown on the first integrals of~\eqref{SIE2a} and~\eqref{SIE2b} means that these integrals are to be evaluated in the Cauchy principle value sense~\cite{PVintegral}.

Equations~\eqref{SIE2a}-\eqref{SIE2b} define the forward scattering problem. For a given $\Gamma=s(x)$ and a given $u^{\mathrm{inc}}(\mathbf{r})$, they are numerically solved for $u(\mathbf{r})$ and $v(\mathbf{r})$, $\mathbf{r} \in \Gamma$ as described next. For numerical solution, the finite domain  $-L/2\leq x\leq L/2$, is divided into $N^{\mathrm{s}}$ number of equal segments of width $w$. Midpoints of these segments are represented by $x_i^{\mathrm{s}}$, $i=1,2,\ldots,N_{\mathrm{s}}$. Unknowns $u(x,s(x))$ and $v(x,s(x))$ are expanded in terms of basis functions as
\begin{subequations}
    \begin{align}
    \label{expansiona} u(x,s(x))&= \sum^{N^{\mathrm{s}}}_{i=1} \bar{u}_i f_i(x)\\
    \label{expansionb} v(x,s(x))&= \sum^{N^{\mathrm{s}}}_{i=1} \bar{v}_i f_i(x).
    \end{align}
    \label{expansion}
\end{subequations}
Here, $f_i(x)$ are the pulse basis functions (see Appendix~\ref{app_a}), and $\bar{u}$ and $\bar{v}$ are the vectors that collect the unknown coefficients associated with these basis functions. Inserting expansions~\eqref{expansiona} and~\eqref{expansionb} into~\eqref{SIE2a} and~\eqref{SIE2b} and point testing the resulting equations at $x_j^{\mathrm{s}}$, $j=1,2,\ldots,N^{\mathrm{s}}$ yields a matrix equation as
\begin{equation}
\label{forward_matrix}\underbrace{\begin{bmatrix}
\bar{Z}^{11} & \bar{Z}^{12} \\
\bar{Z}^{21} & \bar{Z}^{22}
\end{bmatrix}}_{\displaystyle\bar{Z}} \begin{bmatrix} \bar{u}\\ \bar{v}\end{bmatrix} = \begin{bmatrix} \bar{u}^{\mathrm{inc}}\\ \bar{0}\end{bmatrix}.
\end{equation}
Here, $\bar{u}^{\mathrm{inc}}$ is the vector of the tested incident field, and $\bar{Z}$ is the impedance matrix. Their elements are detailed in Appendix~\ref{app_a}. In this work, the matrix equation~\eqref{forward_matrix} is solved by directly inverting the impedance matrix but for large $N_\mathrm{s}$, one can use an iterative method together with well-established acceleration techniques to reduce the computation time and the memory requirement~\cite{FB1,PILE-FBSA,SMCG,FMM_1,FMM_2}.
\subsection{Inverse Scattering Problem}\label{sec:inv_scattering_problem}

The scattered field $u^{\mathrm{sca}}(\mathbf{r})$ in $\Omega_1$ is expressed in terms of $u(\mathbf{r}^\prime)$ and $v(\mathbf{r}^\prime)$ using~\eqref{SIEa}  as
\begin{equation}
\begin{aligned}
     u^{\mathrm{sca}}(\mathbf{r})=u_1(\mathbf{r})-u^{\mathrm{inc}}(\mathbf{r})
    = \int_{\Gamma}\Big[ K_1(\mathbf{r},\mathbf{r}^\prime)u(\mathbf{r}^\prime)-G_1(\mathbf{r},\mathbf{r}^\prime)v(\mathbf{r}^\prime)\Big]dl', \mathbf{r}\in \Omega_1.
\end{aligned}
\label{EqUsSimply}
\end{equation}
Equation~\eqref{EqUsSimply} can be written in a more compact form as
\begin{equation}
    u^{\mathrm{sca}}(\mathbf{r})=\mathcal{D}[s,u,v](\mathbf{r}) \label{operator1}
\end{equation}
where the operator $\mathcal{D}[s,u,v](\mathbf{r})$ is given by
\begin{equation}
    \mathcal{D}[s,u,v](\mathbf{r})=\int_{\Gamma(s)}\Big[ K_1(\mathbf{r},\mathbf{r}^\prime)u(\mathbf{r}^\prime)-G_1(\mathbf{r},\mathbf{r}^\prime)v(\mathbf{r}^\prime)\Big]dl'. \label{operator}
\end{equation}

For the inverse scattering problem, the field scattered from $\Gamma=s(x)$ under excitation by $u^{\mathrm{inc}}(\mathbf{r})$ is measured at points $\mathbf{r}^{\mathrm{r}} = (x_j^{\mathrm{r}},\alpha)$, $j=1,2,\ldots, N^{\mathrm{r}}$ in $\Omega_1$ (see Fig.~\ref{fig:figure1}). This measured scattered field is represented by $u^{\mathrm{mea}}(\mathbf{r}^{\mathrm{r}})$. Then, the inverse scattering problem is defined as reconstructing the unknown $s(x)$ from $u^{\mathrm{mea}}(\mathbf{r}^{\mathrm{r}})$, i.e., it calls for solving
\begin{equation}
\mathcal{D}[s,u,v](x_j^{\mathrm{r}},\alpha)=u^{\mathrm{mea}}(x_j^{\mathrm{r}},\alpha), \,j=1,2,\ldots, N^{\mathrm{r}}
\label{invp}
\end{equation}
for $s(x)$. Equation~\eqref{invp} is nonlinear in $s(x)$, $u(\mathbf{r})$, and $v(\mathbf{r})$ [this can be seen from~\eqref{SIE2a}-\eqref{SIE2b} and~\eqref{operator}-\eqref{operator1}]. Therefore, its numerical solution calls for linearization~\cite{colton}. In this work, this is done using a Newton iterative method~\cite{asTGRS}. Let $n$ represent the Newton iteration number, and superscript ``$(n)$'' attach to a variable in braces mean that that variable is updated/computed at iteration $n$. The resulting Newton update equation reads:
\begin{align}
   \nonumber &\mathcal{D}^\prime\Big[\{s\}^{(n)}, \{u\}^{(n)}, \{v\}^{(n)}\Big](x_j^{\mathrm{r}},\alpha) \{\delta s(x)\}^{(n)} = u^{\mathrm{mea}}(x_j^{\mathrm{r}},\alpha)\\
   \label{NewtonT} &-\mathcal{D}\Big[\{s\}^{(n)}, \{u\}^{(n)}, \{v\}^{(n)}\Big](x_j^{\mathrm{r}},\alpha),\, j=1,2,\ldots, N^{\mathrm{r}}
\end{align}
where $\delta s(x)$ is the unknown Newton step and $\mathcal{D}^\prime[s,u,v](\mathbf{r})$ is the Frechet derivative of the operator $\mathcal{D}[s,u,v](\mathbf{r})$ with respect $s(x)$~\cite{asTGRS}. Equation~\eqref{NewtonT} is numerically solved for $\{\delta s(x)\}^{(n)}$. To facilitate the numerical solution, $\{s(x)\}^{(n)}$ and $\{\delta s(x)\}^{(n)}$ are expanded in terms of (entire-domain) basis functions as 
\begin{subequations}
\label{sp_expansion}\begin{align}
\label{sp_expansiona} \{s(x)\}^{(n)}&=\sum_{i=1}^{N^{\mathrm{p}}} \{\bar{s}_i\}^{(n)}\phi_i(x)\\
\label{sp_expansionb} \{\delta s(x)\}^{(n)}&=\sum_{i=1}^{N^{\mathrm{p}}} \{\bar{d}_i\}^{(n)}\phi_i(x)z
\end{align}
\end{subequations}
Here, $\phi_i(x)$ are the spline-type basis functions (see Appendix~\ref{app_b}), and $\bar{s}$ and $\bar{d}$ are the vectors that collect the coefficients associated with these basis functions. To compute the unknown vector $\bar{d}$, first,~\eqref{sp_expansiona} is used in the forward problem with $\{s(x)\}^{(n)}$ as the input, then the forward problem matrix equation~\eqref{forward_matrix} is solved for $\{\bar{u}\}^{(n)}$ and $\{\bar{v}\}^{(n)}$. $\{u(x,s(x))\}^{(n)}$ and $\{v(x,s(x))\}^{(n)}$ approximated using~\eqref{expansiona} and~\eqref{expansiona}, and $\{s(x)\}^{(n)}$ and $\{\delta s(x)\}^{(n)}$ approximated using~\eqref{sp_expansiona} and~\eqref{sp_expansionb} are inserted into~\eqref{NewtonT}. This yields a matrix equation as
\begin{equation}
\label{inverse_matrix}\{\bar{C}\}^{(n)}\{\bar{d}\}^{(n)} = \{\bar{u}^{\mathrm{mea}}\}^{(n)}-\{\bar{u}^{\mathrm{sca}}\}^{(n)}.
\end{equation}
Here, $\bar{u}^{\mathrm{mea}}$ is the vector that collects the measured scattered field samples, $\bar{u}^{\mathrm{sca}}$ is the vector that collects the samples of the fields scattered from $s(x)$ being reconstructed, and $\bar{C}$ is the matrix that represents the discretized Frechet derivative operator. Their elements are detailed in Appendix~\ref{app_b}.

Inverse scattering problem is ill-posed because measurements are taken at a finite number of points, these measurements are often contaminated by noise, and the integral operator $\mathcal{D}[s,u,v](\mathbf{r})$ has a smoothing effect~\cite{colton,kressReview2018}. This means that the matrix equation~\eqref{inverse_matrix} must be regularized before it can be solved for $\{\bar{d}\}^{(n)}$. To this end, Tikhonov regularization~\cite{asTGRS} is used to convert~\eqref{inverse_matrix} into  
\begin{equation}
\begin{aligned}
\label{reg_inverse_matrix}
(\{\bar{C}^{H}\}^{(n)}\{\bar{C}\}^{(n)}&+\tau\bar{I})\{\bar{d}\}^{(n)}
= \{\bar{C}^{H}\}^{(n)}(\{\bar{u}^{\mathrm{mea}}\}^{(n)}-\{\bar{u}^{\mathrm{sca}}\}^{(n)}).
\end{aligned}
\end{equation}
Here, $\bar{C}^{H}$ is the Hermitian transpose of $\bar{C}$, $\bar{I}$ is the identity matrix, and $\tau$ is the regularization parameter that satisfies $0<\tau<1$. In this work, the matrix equation~\eqref{reg_inverse_matrix} is solved by directly inverting the matrix $(\{\bar{C}^{H}\}^{(n)}\{\bar{C}\}^{(n)}+\tau\bar{I})$. Equation~\eqref{reg_inverse_matrix} is the final discretized form of the Newton update equation that is solved for $\{\bar{d}\}^{(n)}$ at iteration $n$. The next step at iteration $n$ is to update $\{s(x)\}^{(n+1)}$ using
\begin{equation}
\label{update_s}\{s(x)\}^{(n+1)} = \{s(x)\}^{(n)}+\{\delta s(x)\}^{(n)}.
\end{equation}
Equation~\eqref{update_s} can be expressed in terms of the expansion coefficients $\{\bar{s}\}^{(n)}$ and $\{\bar{d}\}^{(n)}$ using~\eqref{sp_expansiona} and~\eqref{sp_expansionb}
\begin{equation}
\label{update_coeff}\{\bar{s}\}^{(n+1)} = \{\bar{s}\}^{(n)}+\{\bar{d}\}^{(n)}.
\end{equation}
Newton iterations are terminated when convergence in reconstruction is achieved, i.e., when the condition
\begin{equation}
\label{convergence} \Big\|\{s(x)\}^{(n+1)}-\{s(x)\}^{(n)}\Big\|_2 = \Big\|\{\delta s(x)\}^{n}\Big\|_2 \leq \xi
\end{equation}
is satisfied. Here, $\xi$ is a user-defined threshold. Condition~\eqref{convergence} can be expressed in terms of the expansion coefficients $\{\bar{s}\}^{(n)}$ and $\{\bar{d}\}^{(n)}$
\begin{equation}
\label{convergence_coeff} \Big\|\{\bar{s}\}^{(n+1)}-\{\bar{s}\}^{(n)}\Big\|_2 = \Big\|\{\bar{d}\}^{n}\Big\|_2 \leq \tilde{\xi},
\end{equation}
where, similarly, $\tilde{\xi}$ is a user-defined threshold. 

\subsection{Multi-Frequency Newton Iterations}\label{sec:multi-freq}
As briefly discussed in Section~\ref{sec:introduction}, multi-frequency data can be included in the iterative reconstruction process described in Section~\ref{sec:inv_scattering_problem} to increase its stability and accuracy~\cite{multiExplain, multiFreqHop, multiTGRSmultiHope}. Assume that the scattered field measurements are taken at frequencies represented by $f_{(m)}$, $m=1,2,\ldots, N^{\mathrm{f}}$. Then, the inverse scattering problem described by~\eqref{invp} is updated as
\begin{align}
\nonumber \mathcal{D}_{(m)}[s,u_{(m)},v_{(m)}](x_j^{\mathrm{r}},\alpha)=u^{\mathrm{mea}}_{(m)}(x_j^{\mathrm{r}},\alpha)\\ 
\label{multi_invp}\,j=1,2,\ldots, N^{\mathrm{r}}, m = 1,2,\ldots, N^{\mathrm{f}}.
\end{align}
In~\eqref{multi_invp} and the rest of the text, subscript ``$(m)$'' attached to a variable means that that variable is updated/computed at frequency $f_{(m)}$. To solve the inverse scattering problem~\eqref{multi_invp}, the Newton method briefly described in Section~\ref{sec:inv_scattering_problem} and detailed in~\cite{asTGRS} is adopted to account for the multi-frequency data. The resulting multi-frequency Newton iterative method reads:
\begin{align*}
&\mathrm{0:\; collect\;} \bar{u}^{\mathrm{mea}}_{(m)}, m = 1,2,\ldots, N^{\mathrm{f}} \mathrm{\;and\;initialize\;} \{\bar{s}_{(1)}\}^{(1)}\\
&\mathrm{1:\; for\;} m=1,2,\ldots, N^{\mathrm{f}} \\
&\mathrm{2.1:\; \;\; for\;} n=1,2,\ldots\\
&\mathrm{2.1.1:\;\; \;\; construct\;} s_{(m)}(x)\mathrm{\;using\;} \{\bar{s}_{(m)}\}^{(n)}\mathrm{\;in}~\eqref{sp_expansiona}\\
&\mathrm{2.1.2:\;\; \;\; discretize\;} s_{(m)}(x) \mathrm{\;using\;} w_{(m)}\\
&\mathrm{2.1.3:\;\; \;\; compute\;} \{\bar{Z}_{(m)}\}^{(n)}, \bar{u}^{\mathrm{inc}}_{(m)}\\
&\mathrm{2.1.4: \;\; \;\; solve}~\eqref{forward_matrix}\mathrm{\;for\;} \{\bar{u}_{(m)}\}^{(n)}, \{\bar{v}_{(m)}\}^{(n)}\\
&\mathrm{2.1.5: \;\; \;\; compute\;}\{\bar{C}_{(m)}\}^{(n)}, \{\bar{u}^{\mathrm{sca}}_{(m)}\}^{(n)} \\
&\mathrm{2.1.6: \;\; \;\; solve}~\eqref{reg_inverse_matrix}\mathrm{\;for\;} \{\bar{d}_{(m)}\}^{(n)}\\
&\mathrm{2.1.7: \;\; \;\; update\;} \{\bar{s}_{(m)}\}^{(n+1)}=\{\bar{s}_{(m)}\}^{(n)}+\{\bar{d}_{(m)}\}^{(n)}\\
&\mathrm{2.1.8: \;\; \;\; check\;}\|\{\bar{d}_{(m)}\}^{(n)}\|_2 \leq \tilde{\xi}\\
&\mathrm{2.1.8.a: \;\; \;\; no:\; continue\;Newton\;iterations\;at\;} f_{(m)} \\
&\mathrm{2.1.8.b: \;\; \;\; yes:\;terminate\;Newton\;iterations\;at\;} f_{(m)}\\
&\mathrm{\;\; \;\; \;\; \;\; \;\; \;\; \;\; \;\;\; update\;} \{\bar{s}_{(m+1)}\}^{(1)}=\{\bar{s}_{(m)}\}^{(n+1)}\\
\end{align*}
Several comments about the above algorithm are in order:
At Step $\mathrm{0}$, $\{\bar{s}_{(1)}\}^{(1)}$ is initialized to zero, i.e., $\{\bar{s}_{(1)}\}^{(1)}=\bar{0}$ which means via~\eqref{sp_expansiona} that the starting point of the iterations is a flat surface at $y=0$. At Step $\mathrm{2.1.2}$, the discretization segment size $w_{(m)}$ is selected based on the spatial variations of the fields at frequency $f_{(m)}$. These are determined by the real and/or imaginary parts of the wavenumbers $k_{1,(m)}$ and $k_{2,(m)}$. Note that as $f_{(m)}$ increases $w_{(m)}$ decreases and the number of segments used in the discretization of the forward scattering problem $N^{\mathrm{s}}_{(m)}=L/w_{(m)}$ increases. At Step $\mathrm{2.1.8}$, the convergence of the Newton iterations $f_{(m)}$ is checked. If the iterations converge (i.e., $\|\{\delta s_{(m)}(x)\}^{(n)}\|_2$ is small enough), the algorithm moves to $f_{(m+1)}$ and sets the initial guess for the Newton iterations at $f_{(m+1)}$ to the reconstruction converged at frequency $f_{(m)}$ (as shown in the second row of Step $\mathrm{2.1.8.b}$).

\section{Numerical Results}\label{sec:numerical_results}
In this section, numerical examples are presented to demonstrate the effectiveness of the multi-frequency Newton iterations. Let $s^{\mathrm{ref}}(x)$ represent the actual surface profile in these numerical examples. $s^{\mathrm{ref}}(x)$ is generated using the stationary random Gaussian process as described in~\cite{thorsos,tsang}. Two user-defined parameters determine the shape of $s^{\mathrm{ref}}(x)$: $\ell$ that denotes the correlation length and $h$ that denotes the standard deviation of the surface height. $s^{\mathrm{ref}}(x)$ is defined in $-L/2\leq x\leq L/2$, where $L=16\,\mathrm{m}$. In addition, tapered cosine windows were implemented on the surface to guarantee local roughness. This is done using the MATLAB  built-in $\mathrm{tukeywin}(L,\, r)$ function, where $r = 25\,\mathrm{cm}$ is the cosine fraction.

To generate $\bar{u}^{\mathrm{mea}}$, first, the fields scattered from $s^{\mathrm{ref}}(x)$ are computed by solving the forward scattering problem, and these fields are synthetically contaminated by noise using
\begin{equation}
    \label{noisyUs}
    \begin{aligned}
    \bar{u}^{\mathrm{mea}}_j=\bar{u}^{\mathrm{ref}}_j+|\bar{u}^{\mathrm{ref}}_j| A_{\mathrm{n}} e^{\mathrm{i}2\pi P_{\mathrm{n}}},\; j = 1,2,\ldots,N^{\mathrm{r}}.
    \end{aligned}
\end{equation}
Here, $A_{\mathrm{n}}$ represents the noise level and $P_{\mathrm{n}}$ is a random number uniformly distributed in the range $[0,\,1]$. 

For all examples considered here, the permittivity and the permeability in $\Omega_1$ and $\Omega_2$ are $\varepsilon_1=\varepsilon_0$ and $\mu_1=\mu_0$ and $\varepsilon_2=4\varepsilon_0$ and $\mu_2=\mu_0$, respectively. Here, $\varepsilon_0$ and $\mu_0$ are the permittivity and the permeability in free space. The conductivity in $\Omega_2$ is $\sigma_2=10^{-5}\,\mathrm{S/m}$. For excitation, the angle of incidence $\theta^{\mathrm{inc}}=0^\circ$ and the taper parameter $g=L/2$. The measurements are taken at points $\mathbf{r}^{\mathrm{r}}=(x_j^{\mathrm{r}}, \alpha)$, $j=1,2, \ldots, N^{\mathrm{r}}$, where $\alpha=4.25\,\mathrm{m}$ and $x^{\mathrm{r}}_j \in [x_{\mathrm{sta}}^{\mathrm{r}}: \Delta x^{\mathrm{r}}: x_{\mathrm{end}}^{\mathrm{r}}]$ with $x_{\mathrm{sta}}^{\mathrm{r}}=-10\,\mathrm{m}$ and $x_{\mathrm{end}}^{\mathrm{r}}=10\,\mathrm{m}$. The samples of frequency of operation $f_{(m)}$, $m=1,2,\ldots,N^{\mathrm{f}}$ are in range $f_{(m)} \in [f_{\mathrm{sta}}: \Delta f: f_{\mathrm{end}}]$.

For all examples, the segment width used in the discretization of the forward scattering problem is $w_{(m)}=2\pi/(10| k_{2,{(m)}}|)$, where $k_{2,{(m)}}$ is the wavenumber in $\Omega_2$ at frequency $f_{(m)}$, the order of the spline basis functions is $p=3$, the Newton iterations are terminated using threshold $\xi=5\times10^{-3}$, and the Tikhonov regularization parameter $\tau=0.75\times10^{-5}$. 

The accuracy of the reconstruction is measured using
\begin{equation}
\label{err_def} {err}_{(m)}=\sqrt{\frac{\displaystyle \sum_{i=1}^{N^{\mathrm{s}}_{(m)}}\Big[s_{(m)}^{(N)}(x^{\mathrm{s}}_{i,(m)})-s^{\mathrm{ref}}(x^{\mathrm{s}}_{i,(m)})\Big]^2}{\displaystyle \sum_{i=1}^{N^{\mathrm{s}}_{(m)}}\Big[s^{\mathrm{ref}}(x^{\mathrm{s}}_{i,(m)})\Big]^2}}
\end{equation}
where $s_{(m)}^{(N)}(x)$ is the reconstruction at the last Newton iteration ($N$) at frequency $f_{(m)}$.  

\subsection{Convergence}
The first example demonstrates the convergence of error in the reconstruction over frequency. For the first example, $s^{\mathrm{ref}}(x)$ is generated with $\ell=0.7\,\mathrm{m}$ and $h=0.07\,\mathrm{m}$. The simulation parameters are $A_{\mathrm{n}}=5\%$,  $f_{\mathrm{sta}}=325\,\mathrm{MHz}$, $f_{\mathrm{end}}=900\,\mathrm{MHz}$, $\Delta f=25\,\mathrm{MHz}$, $\Delta x^{\mathrm{r}}=10\,\mathrm{cm}$, $N^{\mathrm{r}}=200$, and $N^{\mathrm{p}}=25$. 
 
\begin{figure}[ht!]
	\centering
	\includegraphics[width=0.8\columnwidth]{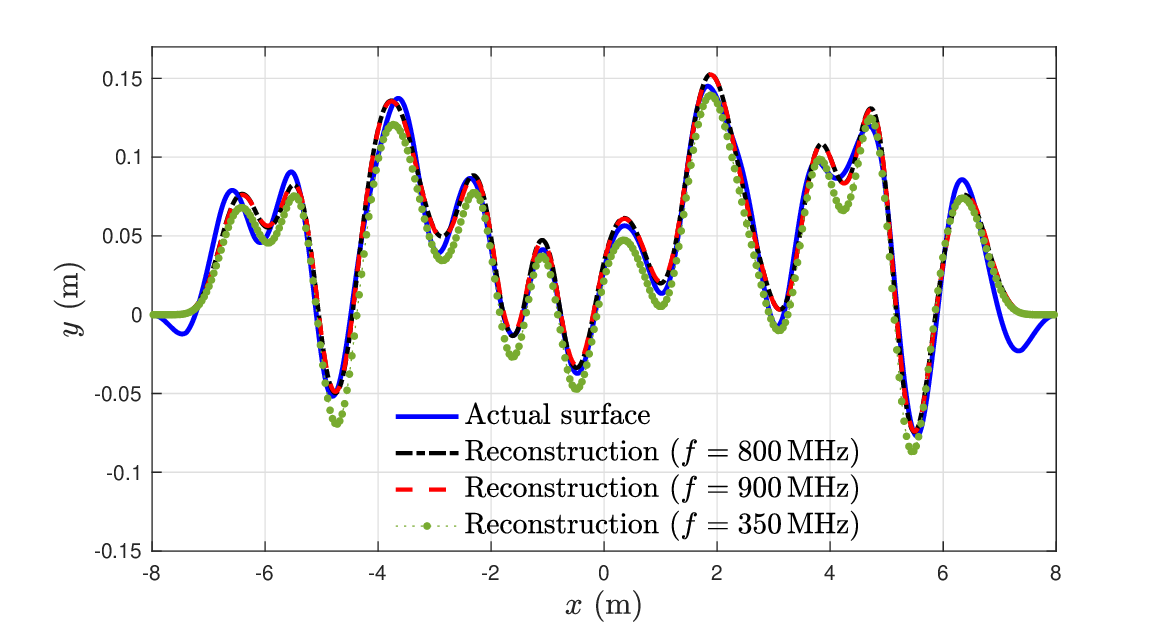}
	\caption{Actual surface profile and the reconstructions at $350\,\mathrm{MHz}$, $800\,\mathrm{MHz}$, and $900\,\mathrm{MHz}$.} \label{fig:figure6}
\end{figure}

\begin{figure}[t!]
	\centering
	\includegraphics[width=0.8\columnwidth]{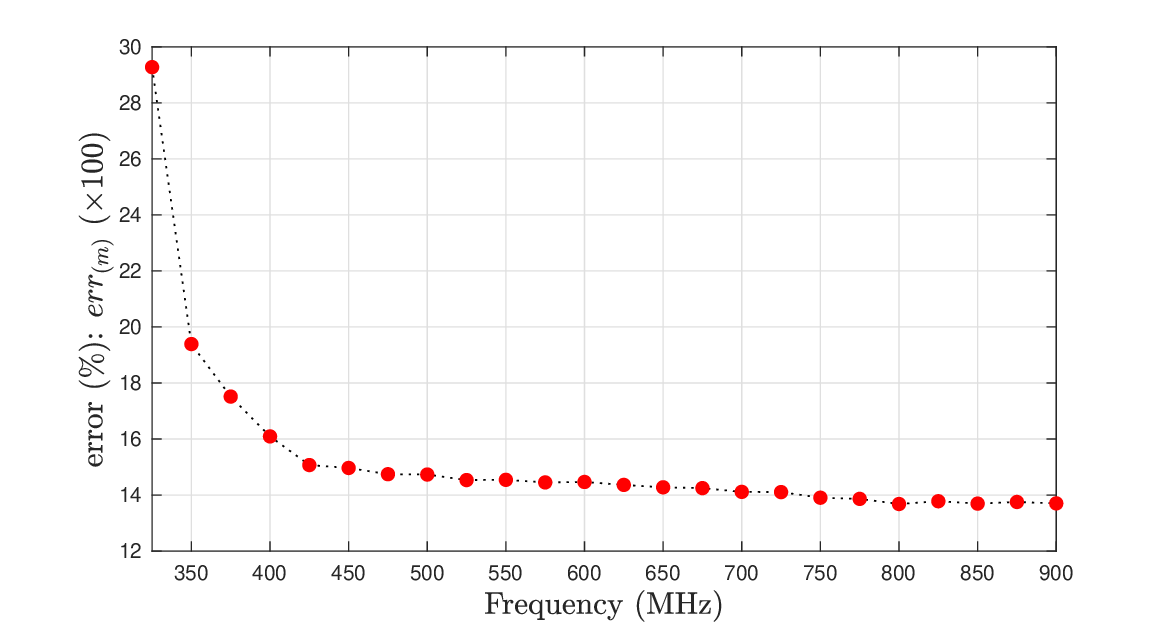}
	\caption{Error in reconstruction computed using~\eqref{err_def} versus frequency.} \label{fig:figure5}
\end{figure}

Fig.~\ref{fig:figure6} compares the surface profiles reconstructed at $350\,\mathrm{MHz}$, $800\,\mathrm{MHz}$, and $900\,\mathrm{MHz}$ to $s^{\mathrm{ref}}(x)$. The figure shows that reconstructions at $800\,\mathrm{MHz}$ and $900\,\mathrm{MHz}$ cannot be easily distinguished by naked eye and are closer to $s^{\mathrm{ref}}(x)$ than the construction at $350\,\mathrm{MHz}$.  Fig.~\ref{fig:figure5} plots ${err}_{(m)}$ versus frequency. It is apparent that the error drops and the quality of the reconstruction
increases with increasing frequency.

 \subsection{Single versus Multi-frequency Reconstruction}
This example demonstrates that the inclusion of multi-frequency data in the reconstruction increases its accuracy. $s^{\mathrm{ref}(x)}$ is generated with $\ell=0.55\,\mathrm{m}$ and $h=0.06\,\mathrm{m}$. Two simulations are carried out. In the first simulation (multi-frequency), the simulation parameters are $A_{\mathrm{n}}=0$ (no noise),  $f_{\mathrm{sta}}=400\,\mathrm{MHz}$, $f_{\mathrm{end}}=600\,\mathrm{MHz}$, $\Delta f=25\,\mathrm{MHz}$, $\Delta x^{\mathrm{r}}=20$, $N^{\mathrm{r}}=100$, and $N^{\mathrm{p}}=20$.  The parameters of the second simulation are same as the first one except the frequency. For the second simulation, the reconstruction is carried out only at a single frequency $600\,\mathrm{MHz}$.

\begin{figure}[t!]
	\centering
	\includegraphics[width=0.8\columnwidth]{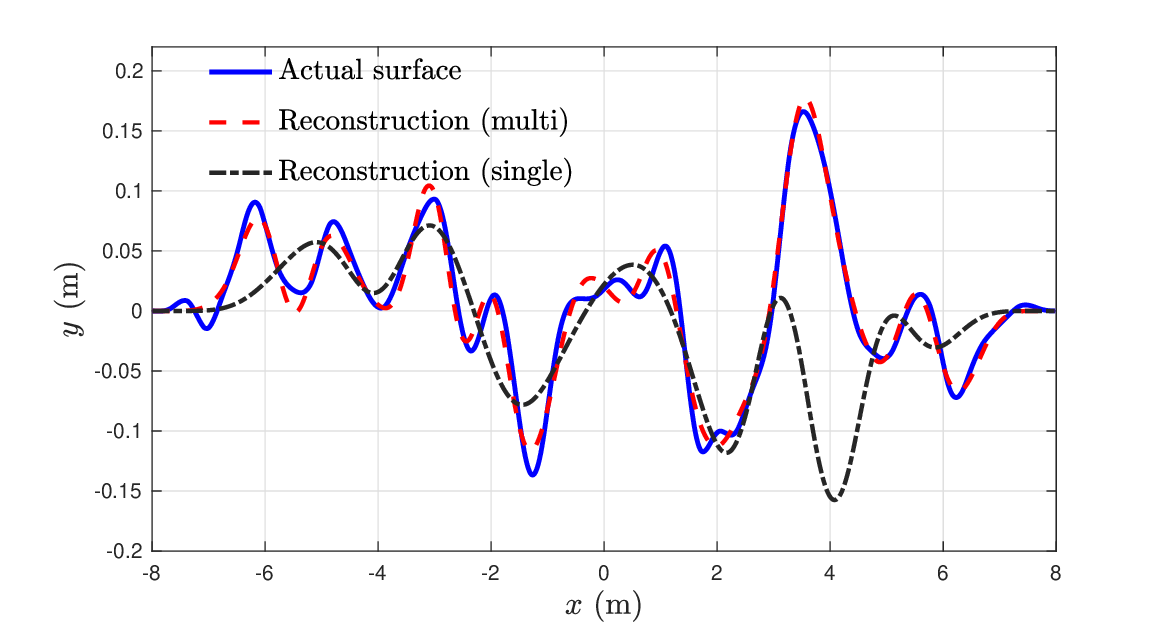}
	\caption{Actual surface profile  and the reconstructions obtained by the multi- and single-frequency simulations at $600\,\mathrm{MHz}$.} \label{fig:figure7}
\end{figure}

Fig.~\ref{fig:figure7} compares the surface profiles reconstructed by the first (multi-frequency) and the second (single-frequency) simulations at $600\,\mathrm{MHz}$ to $s^{\mathrm{ref}}(x)$. As shown in the figure, multi-frequency simulation produces a much more accurate reconstruction. This example clearly demonstrates the benefits of including multi-frequency data in the iterative reconstruction process. 

\subsection{Dependence on Frequency Sampling}
In this example, the dependence of the reconstruction accuracy on the frequency sampling is investigated. $s^{\mathrm{ref}(x)}$ is generated with $\ell=0.4\,\mathrm{m}$ and $h=0.07\,\mathrm{m}$. Five simulations are carried out. The following parameters are kept the same in all five simulations:, $A_{\mathrm{n}}=0$ (no noise),  $f_{\mathrm{sta}}=300\,\mathrm{MHz}$, $f_{\mathrm{end}}=600\,\mathrm{MHz}$, $\Delta x^{\mathrm{r}}=10\,\mathrm{cm}$, $N^{\mathrm{r}}=200$, and $N^{\mathrm{p}}=17$. Then in each simulation $\Delta f$ is set to a different value: $\Delta f \in \{10, 20, 50, 150, 300\}\,\mathrm{MHz}$ (corresponding to $ N^{\mathrm{f}} \in$ $\{31$, $16$, $7$, $3$, $2\}$ number of frequency samples, respectively).

\begin{figure}[t!]
	\centering
	\includegraphics[width=0.8\columnwidth]{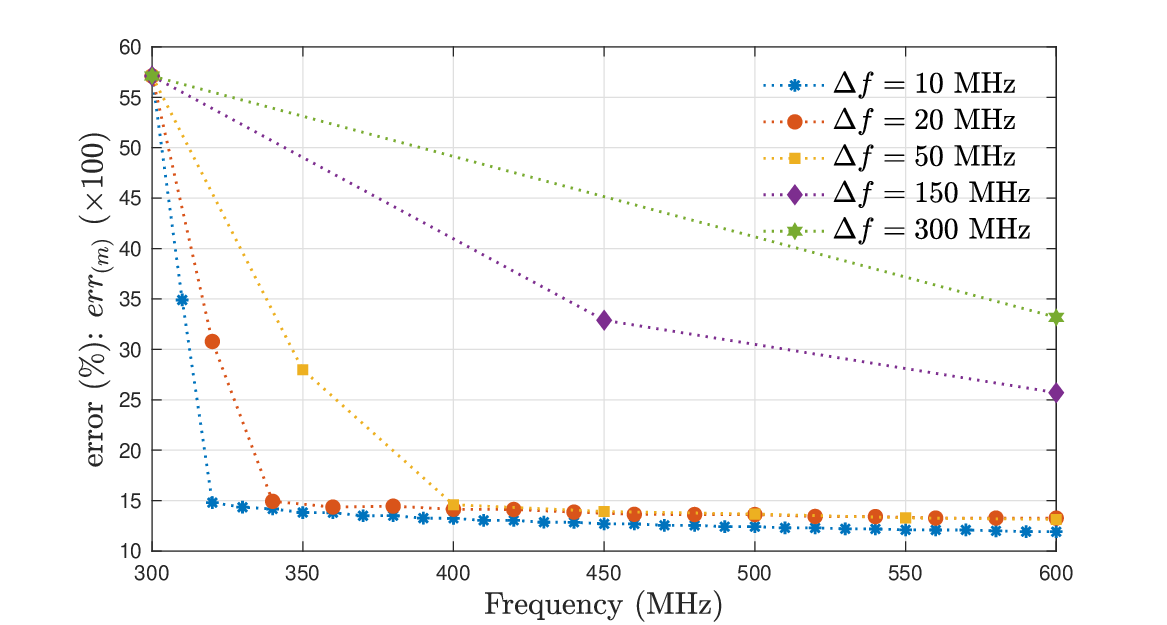}
	\caption{Error in reconstruction computed using~\eqref{err_def} versus frequency for the simulations with the frequency increment $10\,\mathrm{MHz}$, $20\,\mathrm{MHz}$, $50\,\mathrm{MHz}$, $150\,\mathrm{MHz}$, $300\,\mathrm{MHz}$.} \label{fig:figure81}
\end{figure}

\begin{figure}[t!]
	\centering
	\includegraphics[width=0.8\columnwidth]{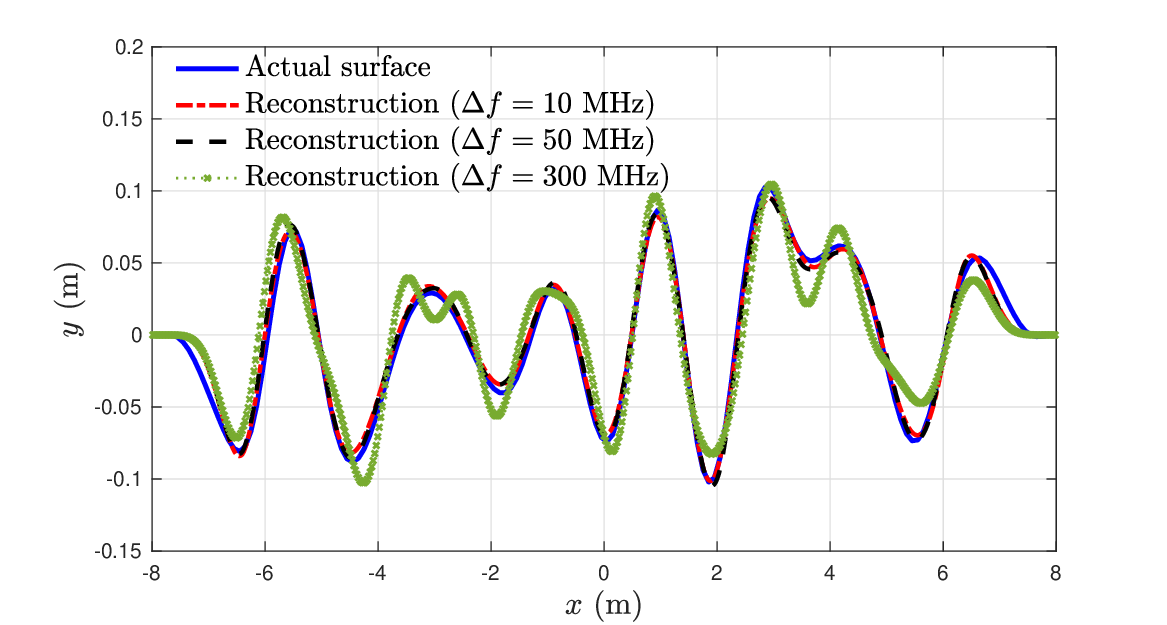}
	\caption{Actual surface profile  and the reconstructions obtained by the simulations with the frequency increment $10\,\mathrm{MHz}$, $50\,\mathrm{MHz}$, and $300\,\mathrm{MHz}$ at $600\,\mathrm{MHz}$.} \label{fig:figure82}
\end{figure}

Fig.~\ref{fig:figure81} plots $err_{(m)}$ versus frequency for all five simulations. The figure shows that smaller $\Delta f$ yields a more accurate reconstruction at $f_{\mathrm{end}}$ (set to $600\,\mathrm{MHz}$ in this case). But this increase in accuracy saturates: The reconstructions at $600\,\mathrm{MHz}$ obtained by simulations with $\Delta f=10\,\mathrm{MHz}$ to $\Delta f=50\,\mathrm{MHz}$ reach roughly the same error level. This is also demonstrated by Fig.~\ref{fig:figure82} where the surface profiles reconstructed by the simulations with $\Delta f=10\,\mathrm{MHz}$, $\Delta f=20\,\mathrm{MHz}$, and $\Delta f=300\,\mathrm{MHz}$ at $600\,\mathrm{MHz}$ are compared to $s^{\mathrm{ref}}(x)$. The reconstructions with $\Delta f=10\,\mathrm{MHz}$ and $\Delta f=50\,\mathrm{MHz}$ are very close to $s^{\mathrm{ref}(x)}$ while the reconstruction with $\Delta f=300\,\mathrm{MHz}$ is not. 

\subsection{Surface Profile with Sharp Variations}
For this example, the performance of the multi-frequency Newton iterations in reconstructing a surface profile with sharp variations. To this end, $s^{\mathrm{ref}}(x)$ is generated using triangular functions 
\begin{equation}
s^{\mathrm{ref}}(x)=
s^{\mathrm{ref}}(x)=
\begin{cases}
  x+4, & \mathrm{if\;} -6 \leq x \leq -3\\
  x,  & \mathrm{if\;}  \quad\;0\leq x \leq 2\\
  4-x,  & \mathrm{if\;} \quad\; 2 \leq x < 4\\
  0, & \mathrm{otherwise}
\end{cases}
\end{equation}
Two simulations are carried out. In the first simulation (multi-frequency), the parameters are $A_{\mathrm{n}}=0$ (no noise),  $f_{\mathrm{sta}}=400\,\mathrm{MHz}$, $f_{\mathrm{end}}=800\,\mathrm{MHz}$, $\Delta f=20\,\mathrm{MHz}$, $\Delta x^{\mathrm{r}}=10\,\mathrm{cm}$, $N^{\mathrm{r}}=200$, and $N^{\mathrm{p}}=18$. The parameters of the second simulation are same as the first one except the frequency. For the second simulation, the reconstruction is carried out only at a single frequency $800\,\mathrm{MHz}$.

\begin{figure}[!t]
	\centering
	\includegraphics[width=0.8\columnwidth]{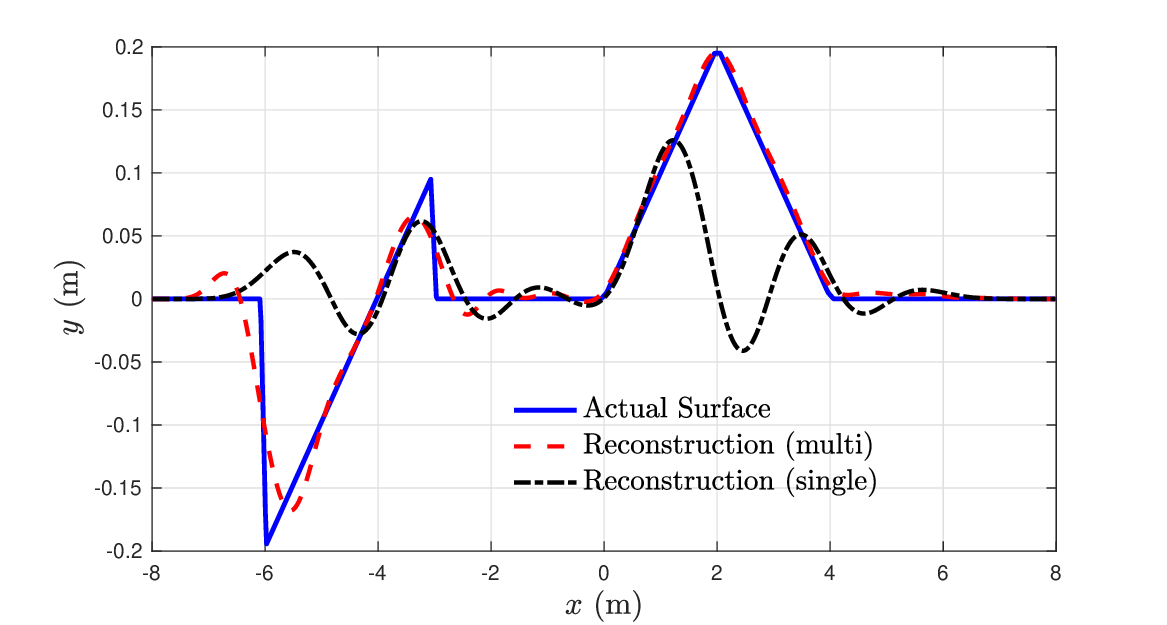}
	\caption{Actual surface profile  and the reconstructions obtained by the multi- and single-frequency simulations at $800\,\mathrm{MHz}$.} \label{fig:triSingleMulti}
\end{figure}

Fig.~\ref{fig:triSingleMulti} compares the surface profiles reconstructed by the first (multi-frequency) and the second (single-frequency) simulations at $800\,\mathrm{MHz}$ to $s^{\mathrm{ref}}(x)$. As expected, multi-frequency reconstruction is significantly more accurate. Single-frequency reconstruction can not capture the sharp variations at all.  

\subsection{Dependence on Noise}
In this example, the performance of the multi-frequency Newton iterations is investigated for measurement data that is contaminated by different levels of noise. $s^{\mathrm{ref}(x)}$ is generated with $\ell=0.5\,\mathrm{m}$ and $h=0.05\,\mathrm{m}$. Eleven simulations are carried out. The following parameters are kept the same in all eleven simulations: $f_{\mathrm{sta}}=425\,\mathrm{MHz}$, $f_{\mathrm{end}}=675\,\mathrm{MHz}$, $\Delta f=25\,\mathrm{MHz}$, $\Delta x^{\mathrm{r}}=10\,\mathrm{cm}$, $N^{\mathrm{r}}=200$, and $N^{\mathrm{p}}=18$. Then in each simulation $A_{\mathrm{n}}$ is set to a different value: $A_{\mathrm{n}} \in$ $\{3$, $5$, $7$, $10$, $12$, $15$, $20$, $25$, $30$, $35$, $40$, $45$, $50\}$$\%$.

\begin{figure}[!t]
	\centering
	\includegraphics[width=0.8\columnwidth]{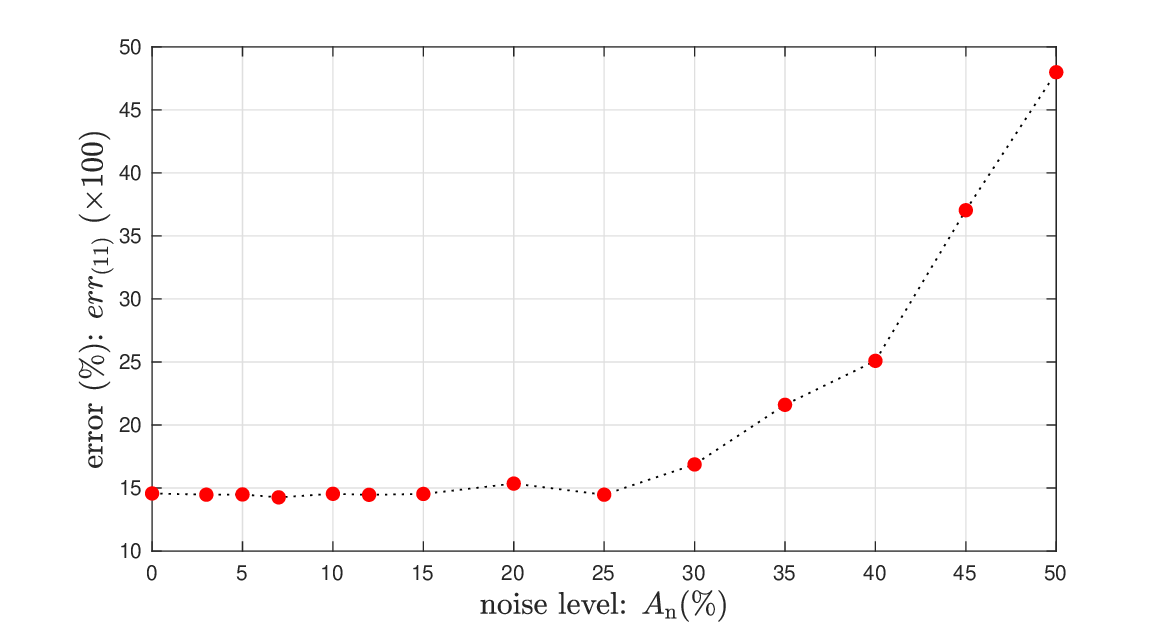}
	\caption{Error in reconstruction at $675\,\mathrm{MHz}$ computed using~\eqref{err_def} versus noise level.} \label{fig:figure8}
\end{figure}

\begin{figure}[!t]
	\centering
	\includegraphics[width=0.8\columnwidth]{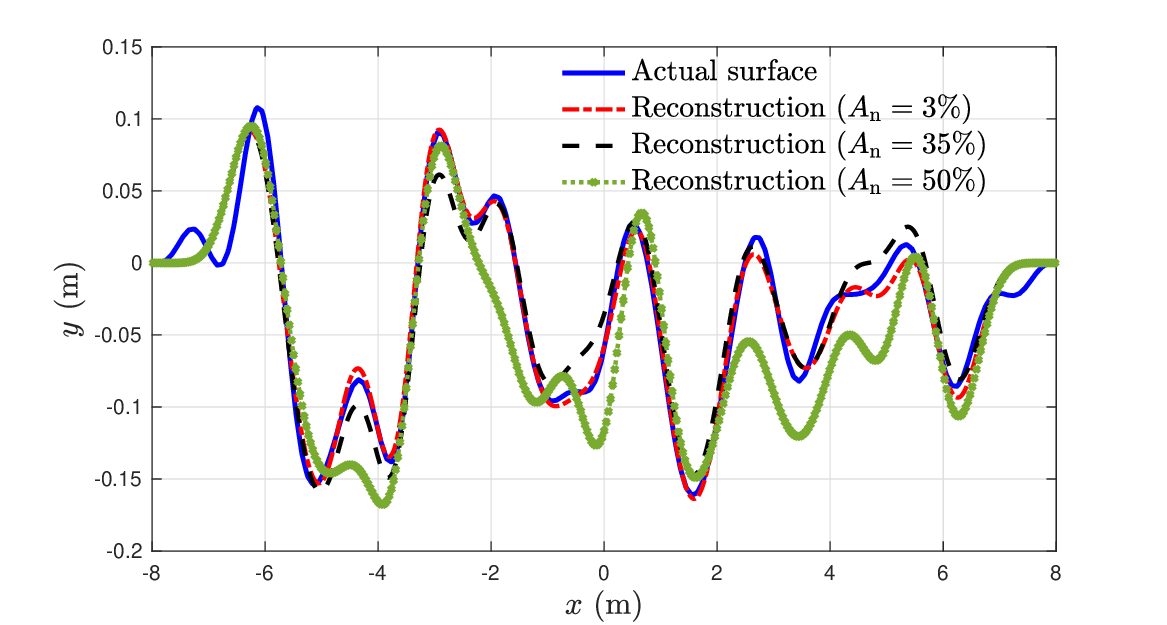}
	\caption{Actual surface profile  and the reconstructions obtained by the simulations with noise level $3\%$, $35\%$, and $50\%$ at $675\,\mathrm{MHz}$.} \label{fig:figure9}
\end{figure}

Fig.~\ref{fig:figure8} plots $err_{11}$ computed at $f_{\mathrm{end}}=675\,\mathrm{MHz}$ by each of these simulations versus the noise level. Fig.~\ref{fig:figure9} compares the surface profiles reconstructed by the simulations with $A_{\mathrm{n}}=3\%$, $A_{\mathrm{n}}=35\%$, and $A_{\mathrm{n}}=50\%$ to $s^{\mathrm{ref}(x)}$. As expected, the measurements that are contaminated the least lead to the most accurate constructions. It can be deducted from Figs.~\ref{fig:figure8} and~\ref{fig:figure9} that the multi-frequency Newton iterations are rather robust to the noise in the measurements.   

\subsection{Dependence on the Number of Measurements}
In the last example, the performance of the multi-frequency Newton iterations is investigated for different numbers of measurement points. $s^{\mathrm{ref}(x)}$ is generated with $\ell=0.7\,\mathrm{m}$ and $h=0.08\,\mathrm{m}$. Twelve simulations are carried out. The following parameters are kept the same in all twelve simulations: $A_{\mathrm{n}}=0$ (no noise),  $f_{\mathrm{sta}}=300\,\mathrm{MHz}$, $f_{\mathrm{end}}=500\,\mathrm{MHz}$, $\Delta f=20\,\mathrm{MHz}$, and $N^{\mathrm{p}}=18$. Then in each simulation $\Delta x^{\mathrm{r}}$ is set to a different value: $\Delta x^{\mathrm{r}} \in $$\{2.5$, $5$, $7.5$, $10$, $15$, $20$, $30$, $40$, $50$, $60$, $80$, $100\}$$\,\mathrm{cm}$ (corresponding to $N^{\mathrm{r}} \in$ $\{800$, $400$, $266$, $200$, $133$, $100$, $66$, $50$, $40$, $33$, $25$, $20\}$ number of measurement points, respectively).

\begin{figure}[t!]
	\centering
	\includegraphics[width=0.8\columnwidth]{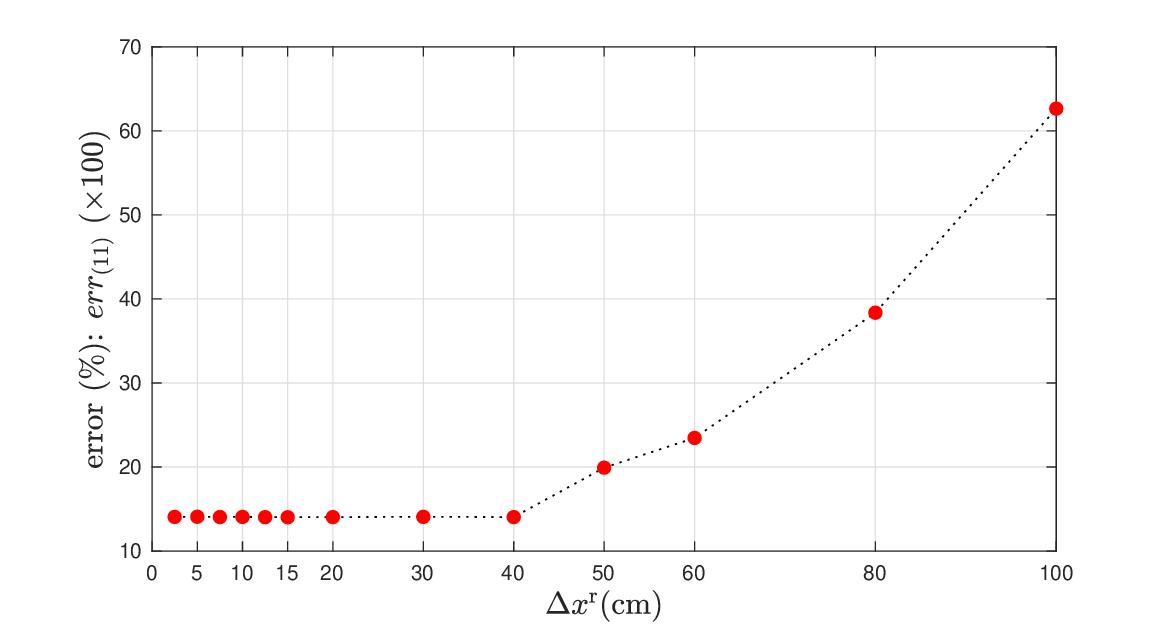}
	\caption{Error in reconstruction at $500\,\mathrm{MHz}$ computed using~\eqref{err_def} versus the measurement point spacing.} \label{fig:figureXX}
\end{figure}

\begin{figure}[t!]
	\centering
	\includegraphics[width=0.8\columnwidth]{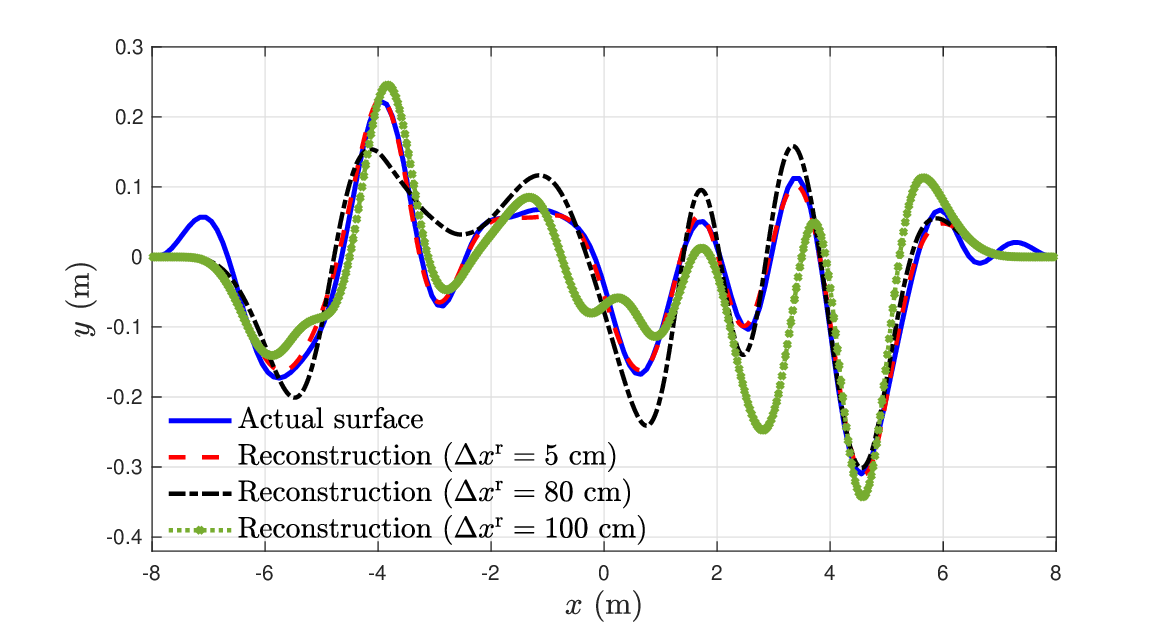}
 \caption{Actual surface profile  and the reconstructions obtained by the simulations with the measurement point spacing $5\,\mathrm{cm}$, $80\,\mathrm{cm}$, and $100\,\mathrm{cm}$ at $500\,\mathrm{MHz}$.}
	\label{fig:figureXX2}
\end{figure}

Fig.~\ref{fig:figureXX} plots $err_{(11)}$ versus frequency for all twelve simulations. Fig.~\ref{fig:figureXX2} where the surface profiles reconstructed by the simulations with $\Delta x^{\mathrm{r}}=5\,\mathrm{cm}$, $\Delta x^{\mathrm{r}}=80\,\mathrm{cm}$, and $\Delta x^{\mathrm{r}}=100\,\mathrm{cm}$ at $500\,\mathrm{MHz}$ are compared to $s^{\mathrm{ref}}(x)$.

\section{Conclusion}\label{sec5conc}
A numerical scheme that utilizes multi-frequency Newton iterations to reconstruct rough surface profiles between two dielectric media is presented. At each frequency sample, the scheme applies Newton iterations to solve the nonlinear inverse scattering problem. In each iteration, the Newton step is determined by solving a linear system that involves the Frechet derivative of the integral operator, which models the scattered fields, and the difference between these fields and the measured data. This linear system is regularized using the Tikhonov method.

The scheme accounts for multi-frequency data in a recursive manner, using the profile reconstructed at one frequency as the initial guess for the next frequency’s iterations. Numerical examples validate the effectiveness of the proposed method, demonstrating its capability to accurately reconstruct surface profiles even in the presence of measurement noise. The results also highlight the superiority of the multi-frequency approach over single-frequency reconstructions, particularly in handling surfaces with sharp variations.

The proposed algorithm can be extended to address layered media problems, where multiple rough surfaces separate more than two media, as well as to 3D surface imaging problems. These extensions are considered for future work.
\setcounter{equation}{0}
\renewcommand{\theequation}{A.\arabic{equation}}
\section*{Appendix A}\label{app_a}
The pulse basis functions $f_i(x)$, $i=1,2,\ldots,N^{\mathrm{s}}$ in expansions~\eqref{expansiona} and~\eqref{expansionb} are defined as
\begin{align}
\label{pulseBasis} f_i(x)= \begin{cases}1 & \mathrm{for\;}(x^{\mathrm{s}}_i-w / 2) \leq x \leq (x^{\mathrm{s}}_i+w/ 2) \\ 0 & \mathrm { otherwise }\end{cases}.
\end{align}
The elements of the vector of tested incident field $\bar{u}^{\mathrm{inc}}$ in~\eqref{forward_matrix} are 
\begin{equation}
\label{uinc} \bar{u}^{\mathrm{inc}}_j = u^{\mathrm{inc}}(x^{\mathrm{s}}_j,s(x^{\mathrm{s}}_j)),\;j = 1,2,\ldots,N^{\mathrm{s}}.
\end{equation}
The midpoint integration is used to evaluate the surface integrals over segment $i$, which arise from the discretization of~\eqref{SIE2a} and~\eqref{SIE2b}. This leads to the following expressions for the elements of the impedance matrix $\bar{Z}$ in~\eqref{forward_matrix}~\cite{tsang}
\begin{subequations}
\begin{align}
\label{z11}&\bar{Z}^{11}_{ji}= 
\begin{cases}
\displaystyle -w_i K_1(\mathbf{r}^{\mathrm{s}}_j,\mathbf{r}^{\mathrm{s}}_i) &\mathrm{\;for\;} j\neq i\\
\displaystyle \frac{1}{2} &\mathrm{\;for\;} j=i
\end{cases}\\
\label{z12}&\bar{Z}^{12}_{ji}=w_i
\begin{cases}
\displaystyle G_1(\mathbf{r}^{\mathrm{s}}_j,\mathbf{r}^{\mathrm{s}}_i) &\mathrm{\;for\;} j\neq i\\
\displaystyle \frac{\mathrm{i}}{4} - \frac{1}{2\pi} \Big[(\gamma - 1) + \ln\Big(\frac{k_1 w_i}{4}\Big)\Big] &\mathrm{\;for\;} j=i
\end{cases}\\
\label{z21}&\bar{Z}^{21}_{ji}=
\begin{cases}
\displaystyle w_i K_2(\mathbf{r}^{\mathrm{s}}_j,\mathbf{r}^{\mathrm{s}}_i) &\mathrm{\;for\;} j\neq i\\
\displaystyle \frac{1}{2} &\mathrm{\;for\;} j=i
\end{cases}\\
\label{z22}&\bar{Z}^{22}_{ji}=-w_i
\begin{cases}
\displaystyle G_2(\mathbf{r}^{\mathrm{s}}_j,\mathbf{r}^{\mathrm{s}}_i) &\mathrm{\;for\;} j\neq i\\
\displaystyle \frac{\mathrm{i}}{4} - \frac{1}{2\pi} \Big[(\gamma - 1) + \ln\Big(\frac{k_2 w_i}{4}\Big)\Big] &\mathrm{\;for\;} j=i
\end{cases}
\end{align}
\end{subequations}
$i,\,j = 1,2,\ldots,N^{\mathrm{s}}$. Here, $w_i=w \sqrt{1+s(x^{\mathrm{s}}_i)^2}$ is the length of the surface corresponding to segment $i$ and $\gamma$ is the Euler's constant. To obtain the expressions for $j=i$ in~\eqref{z12} and~\eqref{z22}, the small argument approximation of the Hankel function $H_0^{(1)}(.)$ is used~\cite{abramowitz}.
\section*{Appendix B}\label{app_b}
The spline-type basis functions $\phi(x)_i$, $i=1,2,\ldots, N^{\mathrm{p}}$ used in expansions~\eqref{sp_expansiona} and~\eqref{sp_expansionb} are defined as~\cite{SplineBook,Li2015,Li2017,zhang2013Spline,asiJRSCNN}
\begin{equation}
    \phi_i(x)=\phi([x-x_i]/g),
\end{equation}
where $g=L/2(N^{\mathrm{p}}+5)$, $x_i=(i+2)(g-L/2)$, and
\begin{equation}
\label{phi_def}\phi(x)= 
\begin{cases}
\displaystyle \begin{aligned}\sum_{q=0}^{p+1} \frac{(-1)^q}{p!}&\binom{p+1}{q}\Big(x+\frac{p+1}{2}-q\Big)^p\\ &\mathrm{for\;}\Big(x+\frac{p+1}{2}\Big) \geq q \end{aligned} \\
\displaystyle 0, \quad\quad\quad\;\,\mathrm{otherwise}
\end{cases}.
\end{equation}	
Here, $p$ is the order of the spline functions.

The elements of the matrix $\bar{C}$ in~\eqref{inverse_matrix} are
\begin{equation}
\label{matrix_Cji}
\begin{aligned}
    C_{ji}=&-\int_{\Gamma(s)}\partial_s G_1(\mathbf{r},\mathbf{r'})\Big|_{\mathbf{r}=(x_j^{\mathrm{r}}, \alpha)}v(\mathbf{r'}) \phi_i(\mathbf{r'})dl'
    +\int_{\Gamma(s)}\partial_s K_1(\mathbf{r},\mathbf{r'})\Big|_{\mathbf{r}=(x_j^{\mathrm{r}}, \alpha)}u(\mathbf{r'}) \phi_i(\mathbf{r'})dl'\\
    &j=1,2,\ldots,N^{\mathrm{r}},\;i=1,2,\ldots,N^{\mathrm{p}}.
\end{aligned}
\end{equation}
Here,``$\partial_s$'' represents the derivative with respective to $s$, and explicit expressions of $\partial_s G_1(\mathbf{r}, \mathbf{r}^{\prime})$ and $\partial_s K_1(\mathbf{r}, \mathbf{r}^{\prime})$ are given in~\cite{asTGRS}. The integral in~\eqref{matrix_Cji} is evaluated using the trapezoidal integration rule. 

The midpoint integration is used to evaluate the surface integrals over segment $i$ which arise from the discretization of $\mathcal{D}[s, u, v](x_j^{\mathrm{r}}, \alpha)$ in~\eqref{NewtonT}. This leads to following expression for the elements of the vector $\bar{u}^{\mathrm{sca}}$ in~\eqref{inverse_matrix}
\begin{equation}
\begin{aligned}
\label{u_bar}\bar{u}^{\mathrm{sca}}_j = \sum_{i=1}^{N^{\mathrm{s}}}\bar{u}_i w_i K_1(\mathbf{r}_j^{\mathrm{r}}, \mathbf{r}_i^{\mathrm{s}})-\sum_{i=1}^{N^{\mathrm{s}}}\bar{v}_i w_i G_1(\mathbf{r}_j^{\mathrm{r}}, \mathbf{r}_i^{\mathrm{s}}),\,
j = 1,2,\ldots, N^{\mathrm{r}}.
\end{aligned}
\end{equation}

\newpage\clearpage

\section*{Figures}

\newpage\clearpage
\bibliographystyle{IEEEtran}
\bibliography{as_ArXiv.bbl}
\end{document}